\documentclass{amsart}
\usepackage{amsthm,amsmath}
\usepackage{amssymb}

\theoremstyle{plain}
\newtheorem{theorem}{Theorem}
\newtheorem{lemma}[theorem]{Lemma}
\newtheorem{proposition}[theorem]{Proposition}
\newtheorem{corollary}[theorem]{Corollary}

\theoremstyle{definition}
\newtheorem{definition}{Definition}
\newtheorem{example}{Example}

\theoremstyle{remark}
\newtheorem{remark}{Remark}

\newcommand\void[1]	{}

\usepackage{amsmath}
\input xypic
\usepackage[all,tips]{xy}
\usepackage{tikz}
\DeclareMathOperator{\Rep}{\mathcal{R}{\sf ep}}
\DeclareMathOperator{\Vect}{\mathcal{V}{\sf ect}}
\DeclareMathOperator{\Aut}{Aut}
\DeclareMathOperator{\ev}{ev}
\DeclareMathOperator{\coev}{coev}
\DeclareMathOperator{\tr}{tr}
\DeclareMathOperator{\Irr}{Irr}

\newcommand{\C}{\mathcal{C}}
\newcommand{\D}{\mathcal{D}}
\newcommand{\E}{\mathcal{E}}
\newcommand{\alg}{\mathrm{alg}}
\newcommand{\loc}{\mathrm{loc}}

\newcommand{\g}{{\mathfrak g}}
\newcommand{\h}{{\mathfrak h}}
\newcommand{\e}{{\mathfrak e}}
\newcommand{\sllie}{\mathfrak{sl}}
\newcommand{\bZ}{\mathbb{Z}}

\newcommand{\tikzmath}[2][.75]
{\vcenter{\hbox{\begin{tikzpicture}[scale=#1] #2\end{tikzpicture}}}
}

\title{Hopf algebra actions in tensor categories}
\author{M.~Bischoff}
\thanks{M.B.\ was supported by NSF DMS Grant 1821162 Quantum Symmetries and Conformal Nets}
\address{ Department of Mathematics\\
Ohio University\\
Athens, OH 45701, USA}
\email{bischoff@ohio.edu}
\author{A.~Davydov}
\address{ Department of Mathematics\\
Ohio University\\
Athens, OH 45701, USA}
\email{davydov@ohio.edu}

\begin{document}

\maketitle

\begin{abstract}
We prove that commutative algebras in braided tensor categories do not admit faithful Hopf algebra actions unless they come from group actions.
We also show that a group action allows us to see the algebra as the regular algebra in the representation category of the acting group. 
\end{abstract}


\section{Introduction}

Hopf algebras are generalisations of group algebras and can be thought of as realisations of ``quantum symmetries" of algebraic objects.
A question of Cohen \cite{co} asks if a commutative algebra can have ``finite quantum symmetries". 
More precisely, the question asks if it is possible for a finite dimensional non-cocommutative semi-simple Hopf algebra to act faithfully on a commutative algebra.
The complete answer is unknown (see however \cite{ew} for a recent progress, which in particular settles in negative the case of inner faithful action).

In the present paper we look into a categorical analog of the Cohen's question.
Namely we examine the ways a Hopf algebra can act faithfully on a separable commutative algebra in a braided tensor category.
We prove that such action could only come from an action by automorphisms. In other words, separable commutative algebras in braided tensor categories do not have interesting quantum symmetries. 

The language of braided tensor categories is proving itself very useful in describing important properties of certain physical systems (e.g. topological orders in condensed matter physics) and goes through the stage of active development.
In particular algebras in braided tensor categories correspond to condensation patterns of a topological order
\cite{ko}.
Other applications of braided tensor categories in quantum field theory are through their relations with conformal nets and vertex operator algebras.

In a recent preprint \cite{dw} Dong and Wang showed that if a finite-dimensional semi-simple Hopf algebra $H$ acts on a vertex operator algebra $V$ (inner) faithfully then the actions comes from a group action. 
In the case when the vertex operator subalgebra of invariants $V^H$ is rational this agrees with our result.
Indeed, according to \cite{hkl} (see also \cite{ckm}) one can see the vertex operator algebra $V$ as an \'etale algebra in the braided tensor category $\Rep(V^H)$ of $V^H$-modules, while the action of $H$ on $V$ translates into an action on that \'etale algebra. 
A similar result has been obtained earlier in the framework of conformal nets in \cite{bi}, which shows that 
finite index \emph{depth two} subnets are given by group fixed points, thus there are no non-trivial faithful 
actions of finite-dimensional C${}^\ast$-Hopf algebras besides the one coming from group algebras. 

We start by reviewing basic facts about separable algebras in braided tensor categories with the emphasis on the their convolution algebras and hypergroups (see \cite{bi} for more details).
Then we define a bialgebra $H$ action on an algebra $A$ in a tensor category and prove that such action gives a homomorphism from the convolution algebra of $A$ to the dual algebra $H^*$. This allows us to show that a faithful bialgebra action on a commutative separable (\'etale) algebra must be a group algebra action.
We conclude by characterising \'etale algebras with a maximal possible automorphism group (maximally symmetric \'etale algebras) in terms of their dimensions. We also show that a maximally symmetric \'etale algebra $A$ in $\C$ gives rise to a braided tensor embedding $F\colon \Rep(G)\to\C$ such that $F$ maps the function algebra $k(G)$ into $A$. Here $G=\Aut_\alg(A)$ is the automorphism group. 

We denote by $k$ a fixed algebraically closed field of characteristic zero.
All our categories will be $k$-linear. We denote the hom-space between objects $X$ and $Y$ of a category $\C$ by $\C(X,Y)$. 
By a tensor category we mean a $k$-linear abelian monoidal category with $k$-linear tensor product.
We denote the monoidal unit object by $I$.  We also assume that the unit object is simple, in particular  $\C(I,I)\cong k$.
By a fusion category we mean a semi-simple spherical tensor category with finitely many (up to isomorphism) simple objects.
We use graphical presentation for morphisms in our braided tensor categories. We read our string diagrams from top to bottom. 

The authors would like to thank Chelsea Walton for useful remarks and the referees for careful reading and helpful suggestions.

\section{\'Etale algebras in braided tensor categories}

Let $\C$ be a spherical tensor category. 
For an object $X\in\C$ denote by $\ev_X\colon X^*\otimes X\to I$ and $\coev_X\colon I\to X\otimes X^*$ the evaluation and the coevaluation morphisms.
Denote by $s_X\colon X\to X^{**}$ the spherical structure morphism.

Let $A=(A,m,\iota)$ be an (associative, unital) algebra in a spherical tensor category $\C$,
where $m\colon A\otimes A\to A$ is the \emph{multiplication}
and $\iota\colon I\to A$ is the \emph{unit} morphisms. 
We call the composite  
$$\xygraph{ !{0;/r6.5pc/:;/u1.5pc/::} 
!{(-1,0)}*{A}="l"  !{(0,0)}*{A\otimes A\otimes A^*}="ml" !{(1,0)}*{A\otimes A^*}="mr" !{(2,0)}*{A^{**}\otimes A^*}="r"  !{(3,0)}*{I}="rr"
"l":"ml"   ^(.4){1\otimes \coev_A}
"ml":"mr" ^{m\otimes1}
"mr":"r"  ^{s_A\otimes 1}
"r":"rr" ^(.6){\mathrm{ev}_{A^*}}
}$$
the {\em canonical trace} of $A$ and denote it by $\varepsilon\colon A\to I$. We call the composite
$$\xymatrix{A \otimes A \ar[rr]^m && A \ar[rr]^{\varepsilon} && I}$$
the {\em canonical pairing} of $A$ and denote it by $b\colon A\otimes A\to I$.
\newline
We call an algebra $A\in\C$ {\em separable} if the canonical pairing is non-degenerate, i.e.\ there is a morphism $\kappa\colon I\to A\otimes A$ such that the composite 
$$\xymatrix{A \ar[rr]^{1\otimes\kappa} && A^{\otimes 3} \ar[rr]^{b\otimes1} && A}$$ 
is the identity. It also implies that the similar composite 
$$\xymatrix{A \ar[rr]^{\kappa\otimes1} && A^{\otimes 3} \ar[rr]^{1\otimes b} && A}$$ 
is also the identity.
\begin{remark}
We prefer this way of defining separability since it make manifest that separability is a property, rather than a structure.
\end{remark}
We use the following graphical representation:
\begin{align*} 
m\ &=\tikzmath{
  \useasboundingbox (-.4,.4) rectangle (1.9,-1.9);
  \draw  (0,0) node [above]  {$\scriptstyle A$} to (0,-.3) to [out=-90,in=-180] (.75,-.8);
\draw  (1.5,0) node [above]  {$\scriptstyle A$} to (1.5,-.3) to [out=-90,in=0] (.75,-.8);
\draw  (.75,-1.5) node [below] {$\scriptstyle A$} to (.75,-.8);
}
&
\iota\ &=
\tikzmath{
  \useasboundingbox(-.4,.9) rectangle (.4,-1.4);
  \draw (0,0)--(0,-1) node [below] {$\scriptstyle A$};
  \node[white] at (0,0)  {$\bullet$};
  \node at (0,0) {$\circ$};
}
\quad.
\end{align*}
and 
$$\varepsilon\ =
\tikzmath{
  \useasboundingbox(-.4,1.4) rectangle (.4,-.9);
  \draw (0,0)--(0,1) node [above] {$\scriptstyle A$};
  \node[white] at (0,0)  {$\bullet$};
  \node at (0,0) {$\circ$};
}
= \qquad
\tikzmath{
\useasboundingbox(-.4,1.4) rectangle (.4,-.9);
\draw (0,1)--(0,.3);
\draw (0,.3) to [out=180,in=90] (-.6,-.3);
\draw (-.6,-.3) to [out=-90,in=180] (0,-.9);
\draw (0,.3) to [out=0,in=90] (.6,-.3);
\draw (.6,-.3) to [out=-90,in=0] (0,-.9);
 }
 \quad.
$$
\\

\begin{remark}
According to our definition $\varepsilon\circ\iota = d(A)1_I$, where $d(A)$ is the dimension of $A\in\C$:
$$
\tikzmath{
  \useasboundingbox(-.4,1) rectangle (.4,-1);
  \draw (0,-.5)--(0,.5) ;
  \node[white] at (0,-.5)  {$\bullet$};
  \node at (0,-.5) {$\circ$};
  \node[white] at (0,.5)  {$\bullet$};
  \node at (0,.5) {$\circ$};
}
\quad = \qquad
\tikzmath{
\useasboundingbox(-.4,1) rectangle (.4,-1);
\draw (0,.6) to [out=180,in=90] (-.6,0);
\draw (-.6,0) to [out=-90,in=180] (0,-.6);
\draw (0,.6) to [out=0,in=90] (.6,0);
\draw (.6,0) to [out=-90,in=0] (0,-.6);
 }
 \quad.
$$
\end{remark}

We call an algebra $A\in\C$ {\em connected} if $\C(I,A) = k$.
\begin{remark}\label{hap}
Any morphism $f:I\to A$ into a connected separable algebra $A$ can be written as $c\iota$, where $c = (\varepsilon\circ f)d(A)^{-1}\in k$\ :
$$
\tikzmath{
  \useasboundingbox(-.4,1) rectangle (.4,-1);
  \draw (0,-.5)--(0,.3) ;
  \draw (-.2,.7) rectangle (.2,.3);
  \node at (0,.5) {$\scriptstyle f$};
}\quad = \quad
\tikzmath{
  \useasboundingbox(-.4,1) rectangle (.4,-1);
  \draw (0,-.5)--(0,.3) ;
  \draw (-.2,.7) rectangle (.2,.3);
  \node at (0,.5) {$\scriptstyle f$};
  \node[white] at (0,-.5)  {$\bullet$};
  \node at (0,-.5) {$\circ$};
}\ 
d(A)^{-1}\ 
\tikzmath{
  \useasboundingbox(-.4,1) rectangle (.4,-1);
  \draw (0,-.5)--(0,.5) ;
  \node[white] at (0,.5)  {$\bullet$};
  \node at (0,.5) {$\circ$};
}
\quad.
$$
\end{remark}

For a separable 
algebra $A\in\C$ we define the {\em convolution algebra}
to be $Q(A)=\C(A,A)$ as a vector space with the multiplication (the convolution product)
$x\ast y = m\circ (x\otimes y)\circ m^\vee$ and the unit
$\iota\circ\varepsilon$. Here $m^\vee\colon A\to A\otimes A$ is the dual morphism to the multiplication.
Graphically 
$$
x\ast y 
\tikzmath{
\node at (-1.3,0) {$=$};
\draw  (.75,1.3) to (.75,.8);
\draw  (0,.2) to (0,.3) to [out=90,in=180] (.75,.8);
\draw (-.2,.2) to (.2,.2) to (.2,.-.2) to (-.2,-.2) to (-.2,.2); 
\node at (0,0) {$\scriptstyle x$};
\draw  (0,-.2) to (0,-.3) to [out=-90,in=-180] (.75,-.8);
\draw   (.75,.8) to [out=0,in=90] (1.5,.3) to (1.5,.2);
\draw (1.3,.2) to (1.7,.2) to (1.7,.-.2) to (1.3,-.2) to (1.3,.2); 
\node at (1.5,0) {$\scriptstyle y$};
\draw  (1.5,-.2) to (1.5,-.3) to [out=-90,in=0] (.75,-.8);
\draw  (.75,-1.3) to (.75,-.8);
}
\quad.
$$

\begin{example}
An algebra endomorphism $g$ of $A$ is an idempotent in the convolution algebra $Q(A)$, i.e. $g*g=g$. 
\end{example}

For an  algebra $A\in\C$ we denote by $\C_A$ the category of its right modules and by ${_A}\C_A$ the category of its bimodules.

Define the map 
\begin{equation}\label{ft}
\phi\colon\C(A,A)\ \to\ {_A}\C_A(A^{\otimes 2},A^{\otimes 2})
\end{equation}
into the space of $A$-bimodule endomorphisms of $A^{\otimes 2}$ by 
$$
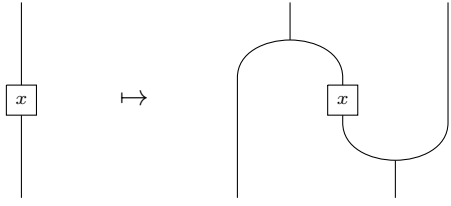
\begin{tikzpicture}
\draw  (0,.2) to (0,1.3);
\draw (-.2,.2) to (.2,.2) to (.2,.-.2) to (-.2,-.2) to (-.2,.2); 
\node at (0,0) {$\scriptstyle x$};
\draw  (0,-.2) to (0,-1.3);
\node at (1.5,0) {\large $\mapsto$};
\end{tikzpicture}
\quad\qquad
\begin{tikzpicture}
\draw  (-.7,1.3) to (-.7,.8);
\draw  (0,.2) to (0,.3) to [out=90,in=0] (-.7,.8);
\draw (-.2,.2) to (.2,.2) to (.2,.-.2) to (-.2,-.2) to (-.2,.2); 
\node at (0,0) {$\scriptstyle x$};
\draw   (-.7,.8) to [out=180,in=90] (-1.4,.3) to (-1.4,-1.3);
\draw  (0,-.2) to (0,-.3) to [out=-90,in=180] (.7,-.8);
\draw (.7,-.8) to [out=0,in=-90] (1.4,-.3) to (1.4,1.3);
\draw  (.7,-1.3) to (.7,-.8);
\end{tikzpicture}
\quad.
$$
We call the map \eqref{ft} the {\em Fourier transform}
\cite{oc}. 

\begin{proposition}
Let $A$ be a separable algebra.
Then the Fourier transform is invertible with the inverse given by
$$
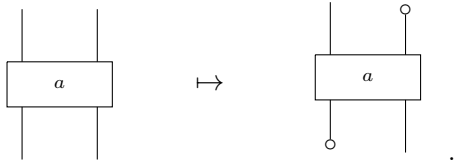
\begin{tikzpicture}
\draw  (-.5,.3) to (-.5,1);
\draw  (.5,.3) to (.5,1);
\draw (-.7,.3) to (.7,.3) to (.7,.-.3) to (-.7,-.3) to (-.7,.3); 
\node at (0,0) {$\scriptstyle a$};
\draw  (.5,.-.3) to (.5,-1);
\draw  (-.5,-.3) to (-.5,-1);
\node at (2,0) {\large $\mapsto$};
\end{tikzpicture}
\quad\qquad
\begin{tikzpicture}
\draw  (-.5,-.3) to (-.5,-.85);
\draw  (.5,.-.3) to (.5,-1);
\draw (-.7,.3) to (.7,.3) to (.7,.-.3) to (-.7,-.3) to (-.7,.3); 
\node at (0,0) {$\scriptstyle a$};
\draw  (.5,.3) to (.5,.85);
\draw  (-.5,.3) to (-.5,1);
\node[white] at (.5,.9) {$\bullet$};
\node at (.5,.9) {$\circ$};
  \node[white] at (-.5,-.9) {$\bullet$};
\node at (-.5,-.9) {$\circ$};
\end{tikzpicture}
\quad.
$$
The Fourier transform has the property $\phi(x*y) = \phi(x)\circ\phi(y)$.
\end{proposition}
\begin{proof}
The invertibility is straightforward.
\newline
The property $\phi(x)\circ\phi(y) = \phi(x*y)$ has the following (also straightforward) graphical verification

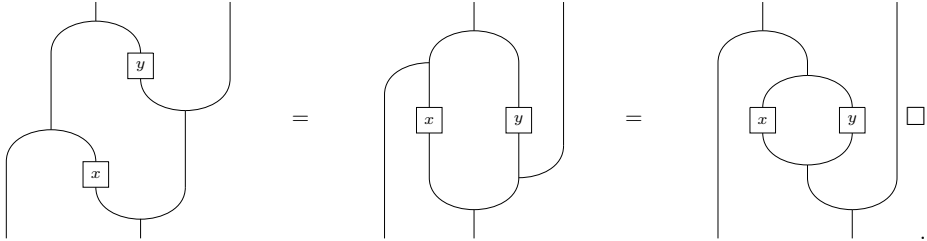
\scalebox{0.85}{\parbox{\linewidth}{%
$$
\begin{tikzpicture}
\draw  (-.7,1) to (-.7,.7);
\draw  (-0,.2) to  [out=90,in=0] (-.7,.7);
\draw (.2,.2) to (-.2,.2) to (-.2,.-.2) to (.2,-.2) to (.2,.2); 
\node at (0,0) {$\scriptstyle y$};
\draw   (-.7,.7) to [out=180,in=90] (-1.4,.2) to (-1.4,-1);
\draw  (0,-.2) to  [out=-90,in=180] (.7,-.7);
\draw (.7,-.7) to [out=0,in=-90] (1.4,-.2) to (1.4,1);
\draw  (.7,-1.9) to (.7,-.7);
\draw (-1.4,-1) to [out=180,in=90] (-2.1,.-1.5) to (-2.1,-2.7);
\draw (-1.4,-1) to [out=0,in=90] (-.7,.-1.5);
\node at (-.7,.-1.7) {$\scriptstyle x$};
\draw (-.5,.-1.5) to (-.9,.-1.5) to (-.9,-1.9) to (-.5,.-1.9) to (-.5,.-1.5); 
\draw  (-.7,-1.9) to  [out=-90,in=180] (0,-2.4) to (0,-2.7);
\draw  (0,-2.4) to  [out=0,in=-90] (.7,-1.9);
\node at (2.5,-.85) {$=$};
\end{tikzpicture}
\qquad\quad
\begin{tikzpicture}
\draw  (1.4,1.85) to (1.4,-.4);
\draw  (-0,1.85) to (-0,1.4);
\draw (1.4,-.4) to [out=-90,in=0] (.7,.-.9);
\draw   (-0,1.4) to [out=0,in=90] (.7,.9) to (.7,.2);
\draw   (-0,1.4) to [out=180,in=90] (-.7,.9) to (-.7,.2);
\draw (.9,.2) to (.5,.2) to (.5,.-.2) to (.9,-.2) to (.9,.2); 
\node at (.7,0) {$\scriptstyle y$};
\draw (-.9,.2) to (-.5,.2) to (-.5,.-.2) to (-.9,-.2) to (-.9,.2); 
\node at (-.7,0) {$\scriptstyle x$};
\draw  (-0,-1.4) to [out=0,in=-90]  (.7,-.9) to (.7,-.2);
\draw  (-0,-1.4) to [out=180,in=-90] (-.7,-.9) to (-.7,-.2);
\draw  (-.7,.9) to [out=180,in=90] (-1.4,.4) to (-1.4,-1.85);
\draw  (-0,-1.85) to (-0,-1.4);
\node at (2.5,0) {$=$};
\end{tikzpicture}
\qquad\quad
\begin{tikzpicture}
\draw  (-.7,1.85) to (-.7,1.4);
\draw   (-.7,1.4) to [out=0,in=90] (-0,.9) to (-0,.7);
\draw   (-0,.7) to [out=0,in=90]  (.7,.2);
\draw   (-0,.7) to [out=180,in=90]  (-.7,.2);
\draw (.9,.2) to (.5,.2) to (.5,.-.2) to (.9,-.2) to (.9,.2); 
\node at (.7,0) {$\scriptstyle y$};
\draw (-.9,.2) to (-.5,.2) to (-.5,.-.2) to (-.9,-.2) to (-.9,.2); 
\node at (-.7,0) {$\scriptstyle x$};
\draw   (-0,-.7) to [out=0,in=-90]  (.7,-.2);
\draw   (-0,-.7) to [out=180,in=-90]  (-.7,-.2);
\draw   (-.7,1.4) to [out=180,in=90] (-1.4,.9) to (-1.4,-1.85);
\draw  (.7,-1.4) to [out=180,in=-90] (-0,-.9) to (-0,-.7);
\draw (.7,-1.4) to [out=0,in=-90] (1.4,-.9) to (1.4,1.85);
\draw  (.7,-1.85) to (.7,-1.4);
\end{tikzpicture}
\quad.
$$
}}
\end{proof}

\begin{corollary}\label{css}
The convolution algebra $Q(A)$ is semi-simple.
\end{corollary}
\begin{proof}
It is known that the category of bimodules over a separable algebra is semi-simple \cite{egno}. 
Now the semi-simplicity of endomorphism algebra ${_A}\C_A(A^{\otimes 2},A^{\otimes 2})$ implies the desired result.
\end{proof}

Let $\C$ be a ribbon tensor category. 
\begin{proposition}\label{prop:Qcomm}
Let $A$ be a commutative separable algebra in a ribbon tensor category $\C$.
Then the convolution algebra $Q(A)$ is commutative.
\end{proposition}
\begin{proof}
Using naturality of the braiding and commutativity of $A$, we get

$$
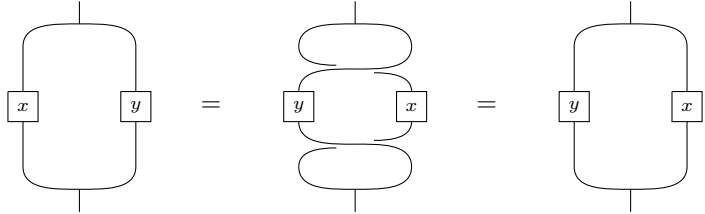
\begin{tikzpicture}
\draw  (.75,1.1) to (.75,1.4);
\draw  (.75,1.1) to  [out=180,in=90] (0,.8) to (0,.2);
\draw  (.75,1.1) to  [out=0,in=90] (1.5,.8) to (1.5,.2);
\draw (-.2,.2) to (.2,.2) to (.2,.-.2) to (-.2,-.2) to (-.2,.2); 
\node at (0,0) {$\scriptstyle x$};
\draw (1.3,.2) to (1.7,.2) to (1.7,.-.2) to (1.3,-.2) to (1.3,.2); 
\node at (1.5,0) {$\scriptstyle y$};
\draw  (1.5,-.2) to (1.5,-.8) to [out=-90,in=0] (.75,-1.1);
\draw  (0,-.2) to (0,-.8) to [out=-90,in=180] (.75,-1.1);
\draw  (.75,-1.4) to (.75,-1.1);
\node at (2.5,0) {$=$};
\end{tikzpicture}
\qquad
\begin{tikzpicture}
\draw  (.75,1.1) to (.75,1.4);
\draw  (.75,1.1) to  [out=180,in=90] (0,.8);
\draw  (0,.8) to  [out=-90,in=180] (.5,.55);
\draw  (1,.45) to  [out=0,in=90] (1.5,.2);
\draw  (.75,1.1) to  [out=0,in=90] (1.5,.8);
\draw  (.75,.5) to  [out=0,in=-90] (1.5,.8);
\draw  (0,.2) to  [out=90,in=180] (.75,.5);
\draw (-.2,.2) to (.2,.2) to (.2,.-.2) to (-.2,-.2) to (-.2,.2); 
\node at (0,0) {$\scriptstyle  y$};
\draw (1.3,.2) to (1.7,.2) to (1.7,.-.2) to (1.3,-.2) to (1.3,.2); 
\node at (1.5,0) {$\scriptstyle  x$};
\draw  (0,-.2) to  [out=-90,in=180] (.75,-.5);
\draw  (.75,-.5) to  [out=0,in=90] (1.5,-.8);
\draw  (.75,-1.1) to  [out=0,in=-90] (1.5,-.8);
\draw  (1,-.45) to  [out=0,in=-90] (1.5,-.2);
\draw  (0,-.8) to  [out=90,in=180] (.5,-.55);
\draw  (.75,-1.1) to  [out=180,in=-90] (0,-.8);
\draw  (.75,-1.1) to (.75,-1.4);
\node at (2.5,0) {$=$};
\end{tikzpicture}
\qquad
\begin{tikzpicture}
\draw  (.75,1.1) to (.75,1.4);
\draw  (.75,1.1) to  [out=180,in=90] (0,.8) to (0,.2);
\draw  (.75,1.1) to  [out=0,in=90] (1.5,.8) to (1.5,.2);
\draw (-.2,.2) to (.2,.2) to (.2,.-.2) to (-.2,-.2) to (-.2,.2); 
\node at (0,0) {$\scriptstyle y$};
\draw (1.3,.2) to (1.7,.2) to (1.7,.-.2) to (1.3,-.2) to (1.3,.2); 
\node at (1.5,0) {$\scriptstyle x$};
\draw  (1.5,-.2) to (1.5,-.8) to [out=-90,in=0] (.75,-1.1);
\draw  (0,-.2) to (0,-.8) to [out=-90,in=180] (.75,-1.1);
\draw  (.75,-1.4) to (.75,-1.1);
\end{tikzpicture}
\quad.
$$
\end{proof}

It follows from Corollary \ref{css} and Proposition \ref{prop:Qcomm} that the convolution algebra $Q(A)$ of a commutative separable algebra $A$ is 
the algebra $k(K)$ of functions on a finite set $K$, the spectrum of $Q(A)$ (which can be defined as the set of homomorphisms $Q(A)\to k$, or equivalently as the set of minimal idempotents).
The composition in $\C(A,A)$ equips the convolution algebra with the second associative multiplication. Its structure constants computed in the basis $K$ 
$$x\circ y = \sum_{z\in K}m^z_{x,y}z\ ,\qquad m^z_{x,y}\in k$$
are invariants of the algebra $A$.
We call the set $K=K(A)$ together with the collection $\{m^z_{x,y}\}_{x,y,z\in K}$ the {\em symmetry hypergroup} of the commutative separable algebra $A$ (see \cite{bi}).

By an {\em \'etale} algebra in $\C$ we mean a commutative, separable algebra such that $\C(I,A)=k$.
In particular, an \'etale algebra is indecomposable.

\begin{proposition}\label{ghg}
Let $A$ be an \'etale algebra and let $g\colon A\to A$ be an algebra automorphism.
\newline
The assignment $x\mapsto tr_A(g\circ x)d(A)^{-1}$ defines an algebra homomorphism $\chi_g\colon Q(A)\to k$.
\newline
Moreover $x*g = \chi_{g^{-1}}(x)g$, so that $g$ is a minimal idempotent in $Q(A)$. 
\end{proposition}
\begin{proof}
Graphically
$$
\begin{tikzpicture}
\node at (-3.5,.35) {$d(A)\chi_g(x) \quad = \quad \tr_A(g\circ x)\quad =$};
\draw  (0,.9) to  [out=90,in=180] (.75,1.4);
\draw (-.2,.9) to (.2,.9) to (.2,.5) to (-.2,.5) to (-.2,.9); 
\node at (0,.7) {$\scriptstyle x$};
\draw  (0,.5) to (0,.2);
\draw (-.2,.2) to (.2,.2) to (.2,.-.2) to (-.2,-.2) to (-.2,.2); 
\node at (0,0) {$\scriptstyle g$};
\draw  (0,-.2) to [out=-90,in=180] (.75,-.7);
\draw   (.75,1.4) to [out=0,in=90] (1.5,.9);
\draw  (1.5,-.2) to  [out=-90,in=0] (.75,-.7);
\draw  (1.5,-.3) to (1.5,.9);
\end{tikzpicture}
\quad.
$$
The homomorphism property $\chi_g(x*y) = \chi_g(x)\chi_g(y)$ has the following graphical justification:
\newline
using homomorphism property of $g$ and commutativity of the multiplication we have
$$
\begin{tikzpicture}
\draw  (.75,2.3) to  [out=0,in=90] (1.5,1.8) to (1.5,-.2);
\draw  (.75,2.3) to  [out=180,in=90] (0,1.8);
\draw  (0,1.8) to  [out=0,in=90] (.7,1.3);
\draw  (0,1.8) to  [out=180,in=90] (-.7,1.3);
\node at (.7,1.1) {$\scriptstyle y$};
\draw (.9,.9) to (.5,.9) to (.5,1.3) to (.9,1.3) to (.9,.9); 
\node at (-.7,1.1) {$\scriptstyle x$};
\draw (-.9,.9) to (-.5,.9) to (-.5,1.3) to (-.9,1.3) to (-.9,.9); 
\draw  (.7,.9) to [out=-90,in=0] (0,.4);
\draw  (-.7,.9) to [out=-90,in=180] (0,.4);
\draw  (0,.4) to (0,.2);
\draw (-.2,.2) to (.2,.2) to (.2,.-.2) to (-.2,-.2) to (-.2,.2); 
\node at (0,0) {$\scriptstyle g$};
\draw  (0,-.2) to [out=-90,in=180] (.75,-.7);
\draw  (1.5,-.2) to  [out=-90,in=0] (.75,-.7);
\node at (2.5,.8) {$=$};
\end{tikzpicture}
\qquad 
\begin{tikzpicture}
\draw  (.75,2.3) to  [out=0,in=90] (1.5,1.8) to (1.5,-.2);
\draw  (.75,2.3) to  [out=180,in=90] (0,1.8);
\draw  (0,1.8) to  [out=0,in=90] (.7,1.3);
\draw  (0,1.8) to  [out=180,in=90] (-.7,1.3);
\node at (.7,1.1) {$\scriptstyle y$};
\draw (.9,.9) to (.5,.9) to (.5,1.3) to (.9,1.3) to (.9,.9); 
\node at (-.7,1.1) {$\scriptstyle x$};
\draw (-.9,.9) to (-.5,.9) to (-.5,1.3) to (-.9,1.3) to (-.9,.9); 
\draw  (.7,.7) to (.7,.9);
\node at (.7,.5) {$\scriptstyle g$};
\draw (.9,.7) to (.5,.7) to (.5,.3) to (.9,.3) to (.9,.7); 
\draw  (-.7,.7) to (-.7,.9);
\node at (-.7,.5) {$\scriptstyle g$};
\draw (-.9,.7) to (-.5,.7) to (-.5,.3) to (-.9,.3) to (-.9,.7); 
\draw  (.7,.3) to [out=-90,in=0] (0,-.2);
\draw  (-.7,.3) to [out=-90,in=180] (0,-.2);
\draw  (0,-.2) to [out=-90,in=180] (.75,-.7);
\draw  (1.5,-.2) to  [out=-90,in=0] (.75,-.7);
\node at (2.5,.8) {$=$};
\end{tikzpicture}
$$
$$
\begin{tikzpicture}
\node at (-2,.8) {$=$};
\draw  (0,1.8) to  [out=0,in=90] (.7,1.3);
\draw  (0,1.8) to  [out=180,in=90] (-.7,1.3);
\node at (.7,1.1) {$\scriptstyle y$};
\draw (.9,.9) to (.5,.9) to (.5,1.3) to (.9,1.3) to (.9,.9); 
\node at (-.7,1.1) {$\scriptstyle x$};
\draw (-.9,.9) to (-.5,.9) to (-.5,1.3) to (-.9,1.3) to (-.9,.9); 
\draw  (.7,.7) to (.7,.9);
\node at (.7,.5) {$\scriptstyle g$};
\draw (.9,.7) to (.5,.7) to (.5,.3) to (.9,.3) to (.9,.7); 
\draw  (-.7,.7) to (-.7,.9);
\node at (-.7,.5) {$\scriptstyle g$};
\draw (-.9,.7) to (-.5,.7) to (-.5,.3) to (-.9,.3) to (-.9,.7); 
\draw  (.7,.3) to [out=-90,in=0] (0,-.2);
\draw  (-.7,.3) to [out=-90,in=180] (0,-.2);
\draw  (0,-.2) to  (0,1.8);
\node at (2,.8) {$=$};
\end{tikzpicture}
\qquad
\begin{tikzpicture}
\draw  (0,1.8) to  [out=0,in=90] (.7,1.3) to (.7,.3);
\draw  (0,1.8) to  [out=180,in=90] (-.7,1.3);
\node at (-.7,1.1) {$\scriptstyle x$};
\draw (-.9,.9) to (-.5,.9) to (-.5,1.3) to (-.9,1.3) to (-.9,.9); 
\draw  (-.7,.7) to (-.7,.9);
\node at (-.7,.5) {$\scriptstyle g$};
\draw (-.9,.7) to (-.5,.7) to (-.5,.3) to (-.9,.3) to (-.9,.7); 
\draw  (.7,.3) to [out=-90,in=0] (0,-.2);
\draw  (-.7,.3) to [out=-90,in=180] (0,-.2);
\draw  (2,1.8) to  [out=0,in=90] (2.7,1.3);
\draw  (2,1.8) to  [out=180,in=90] (1.3,1.3) to (1.3,.3);
\node at (2.7,1.1) {$\scriptstyle y$};
\draw (2.9,.9) to (2.5,.9) to (2.5,1.3) to (2.9,1.3) to (2.9,.9); 
\draw  (2.7,.7) to (2.7,.9);
\node at (2.7,.5) {$\scriptstyle g$};
\draw (2.9,.7) to (2.5,.7) to (2.5,.3) to (2.9,.3) to (2.9,.7); 
\draw  (2.7,.3) to [out=-90,in=0] (2,-.2);
\draw  (1.3,.3) to [out=-90,in=180] (2,-.2);
\draw  (0,-.2) to [out=-90,in=180] (1,-.7);
\draw  (2,-.2) to  [out=-90,in=0] (1,-.7);
\end{tikzpicture}
$$
and the last diagram coincides with (by Remark \ref{hap})
$$
\begin{tikzpicture}
\node at (-1.1,.35) {$d(A)^{-1}$};
\draw  (0,.9) to  [out=90,in=180] (.75,1.4);
\draw (-.2,.9) to (.2,.9) to (.2,.5) to (-.2,.5) to (-.2,.9); 
\node at (0,.7) {$\scriptstyle x$};
\draw  (0,.5) to (0,.2);
\draw (-.2,.2) to (.2,.2) to (.2,.-.2) to (-.2,-.2) to (-.2,.2); 
\node at (0,0) {$\scriptstyle g$};
\draw  (0,-.2) to [out=-90,in=180] (.75,-.7);
\draw   (.75,1.4) to [out=0,in=90] (1.5,.9);
\draw  (1.5,-.2) to  [out=-90,in=0] (.75,-.7);
\draw  (1.5,-.3) to (1.5,.9);
\end{tikzpicture}
\qquad
\begin{tikzpicture}
\draw  (0,.9) to  [out=90,in=180] (.75,1.4);
\draw (-.2,.9) to (.2,.9) to (.2,.5) to (-.2,.5) to (-.2,.9); 
\node at (0,.7) {$\scriptstyle y$};
\draw  (0,.5) to (0,.2);
\draw (-.2,.2) to (.2,.2) to (.2,.-.2) to (-.2,-.2) to (-.2,.2); 
\node at (0,0) {$\scriptstyle g$};
\draw  (0,-.2) to [out=-90,in=180] (.75,-.7);
\draw   (.75,1.4) to [out=0,in=90] (1.5,.9);
\draw  (1.5,-.2) to  [out=-90,in=0] (.75,-.7);
\draw  (1.5,-.3) to (1.5,.9);
\end{tikzpicture}
\quad.
$$
The identity $x*g = \chi_{g^{-1}}(x)g$ is proved as follows
$$
\begin{tikzpicture}
\node at (-2.7,0) {$x*g$}; \node at (-1.3,0) {$=$};
\draw  (.75,1.6) to (.75,.8);
\draw  (0,.2) to (0,.3) to [out=90,in=180] (.75,.8);
\draw (-.2,.2) to (.2,.2) to (.2,.-.2) to (-.2,-.2) to (-.2,.2); 
\node at (0,0) {$\scriptstyle x$};
\draw  (0,-.2) to (0,-.3) to [out=-90,in=-180] (.75,-.8);
\draw   (.75,.8) to [out=0,in=90] (1.5,.3) to (1.5,.2);
\draw (1.3,.2) to (1.7,.2) to (1.7,.-.2) to (1.3,-.2) to (1.3,.2); 
\node at (1.5,0) {$\scriptstyle g$};
\draw  (1.5,-.2) to (1.5,-.3) to [out=-90,in=0] (.75,-.8);
\draw  (.75,-1.6) to (.75,-.8);
\end{tikzpicture}
\qquad
\begin{tikzpicture}
\node at (-1.35,.1) {$=$};
\draw  (.75,1.4) to (.75,1.7);
\draw  (0,.9) to  [out=90,in=180] (.75,1.4);
\draw (-.2,.9) to (.2,.9) to (.2,.5) to (-.2,.5) to (-.2,.9); 
\node at (0,.7) {$\scriptstyle x$};
\draw  (0,.5) to (0,.2);
\draw (-.3,.2) to (.3,.2) to (.3,.-.2) to (-.3,-.2) to (-.3,.2); 
\node at (0,0) {$\scriptstyle  g^{-1}$};
\draw  (0,-.2) to  [out=-90,in=-180] (.75,-.7);
\draw   (.75,1.4) to [out=0,in=90] (1.5,.9);
\draw  (1.5,-.2) to  [out=-90,in=0] (.75,-.7);
\draw  (1.5,-.2) to (1.5,.9);
\draw  (.75,-.7) to (.75,-.9);
\node at (.75,-1.1) {$\scriptstyle  g$};
\draw (.55,-.9) to (.95,-.9) to (.95,-1.3) to (.55,-1.3) to (.55,-.9); 
\draw  (.75,-1.3) to (.75,-1.5);
\end{tikzpicture}
\qquad
$$
$$
\begin{tikzpicture}
\node at (-1.35,-.2) {$=$};
\draw  (2,1.4) to (2,-1.2);
\draw  (0,.9) to  [out=90,in=180] (.75,1.4);
\draw (-.2,.9) to (.2,.9) to (.2,.5) to (-.2,.5) to (-.2,.9); 
\node at (0,.7) {$\scriptstyle x$};
\draw  (0,.5) to (0,.2);
\draw (-.3,.2) to (.3,.2) to (.3,.-.2) to (-.3,-.2) to (-.3,.2); 
\node at (0,0) {$\scriptstyle  g^{-1}$};
\draw  (0,-.2) to  [out=-90,in=-180] (.75,-.7);
\draw   (.75,1.4) to [out=0,in=90] (1.5,.9);
\draw  (1.5,-.2) to  [out=-90,in=0] (.75,-.7);
\draw  (1.5,-.2) to (1.5,.9);
\draw  (.75,-.7) to  [out=-90,in=180] (2,-1);
\node at (2,-1.4) {$\scriptstyle  g$};
\draw (1.8,-1.2) to (2.2,-1.2) to (2.2,-1.6) to (1.8,-1.6) to (1.8,-1.2); 
\draw  (2,-1.6) to (2,-1.8);
\end{tikzpicture}
\qquad
\begin{tikzpicture}
\node at (-1.6,-.2) {$=\quad d(A)^{-1}$};
\draw  (2,1.4) to (2,-1.2);
\draw  (0,.9) to  [out=90,in=180] (.75,1.4);
\draw (-.2,.9) to (.2,.9) to (.2,.5) to (-.2,.5) to (-.2,.9); 
\node at (0,.7) {$\scriptstyle x$};
\draw  (0,.5) to (0,.2);
\draw (-.3,.2) to (.3,.2) to (.3,.-.2) to (-.3,-.2) to (-.3,.2); 
\node at (0,0) {$\scriptstyle  g^{-1}$};
\draw  (0,-.2) to  [out=-90,in=-180] (.75,-.7);
\draw   (.75,1.4) to [out=0,in=90] (1.5,.9);
\draw  (1.5,-.2) to  [out=-90,in=0] (.75,-.7);
\draw  (1.5,-.2) to (1.5,.9);
\node at (2,-1.4) {$\scriptstyle  g$};
\draw (1.8,-1.2) to (2.2,-1.2) to (2.2,-1.6) to (1.8,-1.6) to (1.8,-1.2); 
\draw  (2,-1.6) to (2,-1.8);
\end{tikzpicture}
\qquad.
$$

\void{
$$
\tikzmath{
\draw  (.75,1.6) to (.75,.8);
\draw  (0,.2) to (0,.3) to [out=90,in=180] (.75,.8);
\draw (-.2,.2) to (.2,.2) to (.2,.-.2) to (-.2,-.2) to (-.2,.2); 
\node at (0,0) {$\scriptstyle x$};
\draw  (0,-.2) to (0,-.3) to [out=-90,in=-180] (.75,-.8);
\draw   (.75,.8) to [out=0,in=90] (1.5,.3) to (1.5,.2);
\draw (1.3,.2) to (1.7,.2) to (1.7,.-.2) to (1.3,-.2) to (1.3,.2); 
\node at (1.5,0) {$\scriptstyle g$};
\draw  (1.5,-.2) to (1.5,-.3) to [out=-90,in=0] (.75,-.8);
\draw  (.75,-1.6) to (.75,-.8);
}
~=~
\tikzmath{
\draw  (.75,1.4) to (.75,1.7);
\draw  (0,.9) to  [out=90,in=180] (.75,1.4);
\draw (-.2,.9) to (.2,.9) to (.2,.5) to (-.2,.5) to (-.2,.9); 
\node at (0,.7) {$\scriptstyle x$};
\draw  (0,.5) to (0,.2);
\draw (-.3,.2) to (.3,.2) to (.3,.-.2) to (-.3,-.2) to (-.3,.2); 
\node at (0,0) {$\scriptstyle  g^{-1}$};
\draw  (0,-.2) to  [out=-90,in=-180] (.75,-.7);
\draw   (.75,1.4) to [out=0,in=90] (1.5,.9);
\draw  (1.5,-.2) to  [out=-90,in=0] (.75,-.7);
\draw  (1.5,-.2) to (1.5,.9);
\draw  (.75,-.7) to (.75,-.9);
\node at (.75,-1.1) {$\scriptstyle  g$};
\draw (.55,-.9) to (.95,-.9) to (.95,-1.3) to (.55,-1.3) to (.55,-.9); 
\draw  (.75,-1.3) to (.75,-1.5);
}
~=~
\tikzmath{
\draw  (2,1.4) to (2,-1.2);
\draw  (0,.9) to  [out=90,in=180] (.75,1.4);
\draw (-.2,.9) to (.2,.9) to (.2,.5) to (-.2,.5) to (-.2,.9); 
\node at (0,.7) {$\scriptstyle x$};
\draw  (0,.5) to (0,.2);
\draw (-.3,.2) to (.3,.2) to (.3,.-.2) to (-.3,-.2) to (-.3,.2); 
\node at (0,0) {$\scriptstyle  g^{-1}$};
\draw  (0,-.2) to  [out=-90,in=-180] (.75,-.7);
\draw   (.75,1.4) to [out=0,in=90] (1.5,.9);
\draw  (1.5,-.2) to  [out=-90,in=0] (.75,-.7);
\draw  (1.5,-.2) to (1.5,.9);
\draw  (.75,-.7) to  [out=-90,in=180] (2,-1);
\node at (2,-1.4) {$\scriptstyle  g$};
\draw (1.8,-1.2) to (2.2,-1.2) to (2.2,-1.6) to (1.8,-1.6) to (1.8,-1.2); 
\draw  (2,-1.6) to (2,-1.8);
}
~=~\quad d(A)^{-1}\cdot
\tikzmath{
\draw  (2,1.4) to (2,-1.2);
\draw  (0,.9) to  [out=90,in=180] (.75,1.4);
\draw (-.2,.9) to (.2,.9) to (.2,.5) to (-.2,.5) to (-.2,.9); 
\node at (0,.7) {$\scriptstyle x$};
\draw  (0,.5) to (0,.2);
\draw (-.3,.2) to (.3,.2) to (.3,.-.2) to (-.3,-.2) to (-.3,.2); 
\node at (0,0) {$\scriptstyle  g^{-1}$};
\draw  (0,-.2) to  [out=-90,in=-180] (.75,-.7);
\draw   (.75,1.4) to [out=0,in=90] (1.5,.9);
\draw  (1.5,-.2) to  [out=-90,in=0] (.75,-.7);
\draw  (1.5,-.2) to (1.5,.9);
\node at (2,-1.4) {$\scriptstyle  g$};
\draw (1.8,-1.2) to (2.2,-1.2) to (2.2,-1.6) to (1.8,-1.6) to (1.8,-1.2); 
\draw  (2,-1.6) to (2,-1.8);
}
\quad.
$$
}
\end{proof}

Proposition \ref{ghg} says that the automorphism group $\Aut_{\alg}(A)$ is a subset of its symmetry hypergroup $K(A)$ and that the structure constants of $\Aut_{\alg}(A)$ are given by the group operation, i.e.\ that  the automorphism group $\Aut_{\alg}(A)$ is a sub-hypergroup of $K(A)$.
\newline
We call an \'etale algebra $A$ {\em Galois} if $\Aut_{\alg}(A) = K(A)$, i.e.\ if $\C(A,A) = k[\Aut_{\alg}(A)]$.
In Section \ref{gal} we give a convenient criterion for being Galois.

\section{Bialgebra actions on \'etale algebras}

Let $A$ be an algebra in a tensor category $\C$. Let $H$ be a bialgebra in the category of vector spaces $\Vect$, which we consider as the tensor subcategory of (the symmetric centre of the monoidal centre of) $\C$. 
An {\em action} of a bialgebra $H$ on $A$ is a morphism $a\colon H\otimes A\to A$ such that the diagrams
$$\xymatrix{H\otimes H\otimes A \ar[r]^{a\otimes 1} \ar[d]_{m_H\otimes 1}& H\otimes A \ar[d]^a \\ H\otimes A \ar[r]^a & A}\qquad\qquad
\xymatrix{ A \ar@{=}[rd] \ar[d]_{\iota_H1} \\ H\otimes A \ar[r]^a & A
}$$
$$\xymatrix{
H\otimes A\otimes A \ar[r]^{1\otimes m_A} \ar[d]_{\delta\otimes 1} & H\otimes A\ar[r]^a & A\\
H\otimes H\otimes A\otimes A \ar[r]^{1\otimes c\otimes 1} & H\otimes A\otimes H\otimes A \ar[r]^(.6){a\otimes a} & A\otimes A \ar[u]_{m_A}}\qquad\qquad
\xymatrix{ H \ar[r]^{1\iota} \ar[d]_\varepsilon & H\otimes A \ar[d]^a\\ I \ar[r]^\iota & A
}$$
commute. 
Here $\iota_H\colon k\to H$ is the unit of $H$, $c\colon H\otimes A\to A\otimes H$ is the braiding of a vector space $H$ with an object $A\in\C$, and $\ \delta\colon H\to H\otimes H$ is the coproduct of $H$. 
\begin{remark}
Note that the first diagram says that $A$ is an $H$-module.
The second diagram makes this $H$-module unital.
The last two diagrams say that the multiplication and the unit morphisms of $A$ are homomorphisms of $H$-modules.
\end{remark}
Graphically the third condition has the form
$$
\begin{tikzpicture}
\draw  (0,.4) to (0,.6); 
\draw (.5,-.4) to (.9,-.4) to (.9,.-.8) to (.5,-.8) to (.5,-.4); 
\node at (.7,-.6) {$\scriptstyle h$};
\draw  (1.4,.4) to (1.4,.6); 
\draw  (0,.4) to [out=-90,in=180] (.7,-.1);
\draw (.7,-.1) to [out=0,in=-90] (1.4,.4) ;
\draw  (.7,-1.3) to (.7,-.8); 
\draw  (.7,-.4) to (.7,-.1);
\node at (2.4,-.3) {\large${ = }$};
\end{tikzpicture}
\qquad
\begin{tikzpicture}
\node at (-1,-.3) {\large ${ \sum_{(h)}}$};
\draw  (0,.2) to (0,.6); 
\draw (-.3,.2) to (.3,.2) to (.3,.-.2) to (-.3,-.2) to (-.3,.2); 
\node at (0,0) {$\scriptstyle h_{(0)}$};
\draw (1.1,.2) to (1.7,.2) to (1.7,.-.2) to (1.1,-.2) to (1.1,.2); 
\node at (1.4,0) {$\scriptstyle h_{(1)}$};
\draw  (1.4,.2) to (1.4,.6); 
\draw  (0,-.2) to (0,-.3) to [out=-90,in=180] (.7,-.8);
\draw (.7,-.8) to [out=0,in=-90] (1.4,-.2) ;
\draw  (.7,-1.3) to (.7,-.8);
\end{tikzpicture}
\quad.
$$

\void{
$$
\tikzmath{
\draw  (0,.4) to (0,.6); 
\draw (.5,-.4) to (.9,-.4) to (.9,.-.8) to (.5,-.8) to (.5,-.4); 
\node at (.7,-.6) {$\scriptstyle h$};
\draw  (1.4,.4) to (1.4,.6); 
\draw  (0,.4) to [out=-90,in=180] (.7,-.1);
\draw (.7,-.1) to [out=0,in=-90] (1.4,.4) ;
\draw  (.7,-1.3) to (.7,-.8); 
\draw  (.7,-.4) to (.7,-.1);
}
~=~
 \sum\limits_{(h)}
\tikzmath{
\draw  (0,.2) to (0,.6); 
\draw (-.3,.2) to (.3,.2) to (.3,.-.2) to (-.3,-.2) to (-.3,.2); 
\node at (0,0) {$\scriptstyle h_{(0)}$};
\draw (1.1,.2) to (1.7,.2) to (1.7,.-.2) to (1.1,-.2) to (1.1,.2); 
\node at (1.4,0) {$\scriptstyle h_{(1)}$};
\draw  (1.4,.2) to (1.4,.6); 
\draw  (0,-.2) to (0,-.3) to [out=-90,in=180] (.7,-.8);
\draw (.7,-.8) to [out=0,in=-90] (1.4,-.2) ;
\draw  (.7,-1.3) to (.7,-.8);
}
\quad.
$$
}
Here we use Sweedler's notation for the comultiplication $\delta(h) = \sum_{(h)}h_{(0)}\otimes h_{(1)}$. 
\newline
Note that an action of $H$ can be rewritten as a linear map $H\to \C(A,A)$, which in particular is a homomorphism of algebras (with respect to the composition on $\C(A,A)$).
\newline
We say that an action is {\em faithful} if the corresponding map $H\to \C(A,A)$ is an embedding.

\begin{example}
Let $G\subset \Aut_{\alg}(A)$ be a subgroup. Then by linear extension we get a  Hopf action of
the group algebra $k[G]$ on $A$. 
\end{example}

Denote by $H^*$ the dual Hopf algebra of $H$.
The multiplication on $H^\ast$ is given 
by 
$$(l\cdot m)(h)=\sum_{(h)}l(h_{(0)})m(h_{(1)})  \quad\quad h\in H\,,\quad  l,m\in H^*\,.$$

\begin{proposition}
\label{prop:Qquotient}
Let $A\in\C$ be a separable connected algebra, and $H$ a  Hopf algebra faithfully 
acting on $A$. Then there is an epimorphism $\gamma\colon Q(A)\to H^\ast$. 
In particular, $H^\ast$ is a quotient of $Q(A)$. 
\end{proposition}
\begin{proof}
Define the pairing $\eta\colon\C(A,A)\times \C(A,A) \to k$ by $\eta(a,b)=tr_A(b\circ a)$. 
Note that this pairing is defined for any object $A\in \C$ and that it is the direct sum of pairings for isotypical components of $A$.
On each isotypical component, i.e.\ on a direct sum of a simple object $X$ in $\C$, the pairing is proportional to the canonical pairing on a matrix algebra (with the proportionality coefficient being $d(X)$).
Thus the pairing $\eta$ is non-degenerate.
\newline
Graphically
$$
\begin{tikzpicture}
\node at (-2.7,.35) {$\tr_A(b\circ a)$};
\node at (-1.3,.35) {$=$};
\draw  (0,.9) to  [out=90,in=180] (.75,1.4);
\draw (-.2,.9) to (.2,.9) to (.2,.5) to (-.2,.5) to (-.2,.9); 
\node at (0,.7) {$\scriptstyle a$};
\draw  (0,.5) to (0,.2);
\draw (-.2,.2) to (.2,.2) to (.2,.-.2) to (-.2,-.2) to (-.2,.2); 
\node at (0,0) {$\scriptstyle b$};
\draw  (0,-.2) to  [out=-90,in=-180] (.75,-.7);
\draw   (.75,1.4) to [out=0,in=90] (1.5,.9);
\draw  (1.5,-.2) to [out=-90,in=0] (.75,-.7);
\draw  (1.5,-.2) to (1.5,.9);
\end{tikzpicture}
\quad.
$$
Define a surjective $k$-linear map $\gamma\colon Q(A)\to H^\ast$  by 
$$\gamma(a)(h) = d(A)^{-1}\eta(a,h),\qquad a\in Q(A),\ h\in H\ .$$
By the definition of the  multiplication on $H^\ast$, the homomorphism property 
\newline
$\gamma(a*b)=\gamma(a)\cdot\gamma(b)$ is equivalent to 
$$\eta(a*b,h) = d(A)^{-1}\sum_{(h)}\eta(a,h_{(0)})\eta(b,h_{(1)})\ .$$
The last identity has the following graphical verification:
\newline
\scalebox{0.85}{\parbox{\linewidth}{%
$$
\begin{tikzpicture}
\draw  (.75,2.3) to  [out=0,in=90] (1.5,1.8) to (1.5,-.2);
\draw  (.75,2.3) to  [out=180,in=90] (0,1.8);
\draw  (0,1.8) to  [out=0,in=90] (.7,1.3);
\draw  (0,1.8) to  [out=180,in=90] (-.7,1.3);
\node at (.7,1.1) {$\scriptstyle b$};
\draw (.9,.9) to (.5,.9) to (.5,1.3) to (.9,1.3) to (.9,.9); 
\node at (-.7,1.1) {$\scriptstyle a$};
\draw (-.9,.9) to (-.5,.9) to (-.5,1.3) to (-.9,1.3) to (-.9,.9); 
\draw  (.7,.9) to [out=-90,in=0] (0,.4);
\draw  (-.7,.9) to [out=-90,in=180] (0,.4);
\draw  (0,.4) to (0,.2);
\draw (-.2,.2) to (.2,.2) to (.2,.-.2) to (-.2,-.2) to (-.2,.2); 
\node at (0,0) {$\scriptstyle h$};
\draw  (0,-.2) to [out=-90,in=180] (.75,-.7);
\draw  (1.5,-.2) to  [out=-90,in=0] (.75,-.7);
\node at (3,.8) {$=\quad\ \sum_{(h)}$};
\end{tikzpicture}
\quad 
\begin{tikzpicture}
\draw  (.75,2.3) to  [out=0,in=90] (1.5,1.8) to (1.5,-.2);
\draw  (.75,2.3) to  [out=180,in=90] (0,1.8);
\draw  (0,1.8) to  [out=0,in=90] (.7,1.3);
\draw  (0,1.8) to  [out=180,in=90] (-.7,1.3);
\node at (.7,1.1) {$\scriptstyle b$};
\draw (.9,.9) to (.5,.9) to (.5,1.3) to (.9,1.3) to (.9,.9); 
\node at (-.7,1.1) {$\scriptstyle a$};
\draw (-.9,.9) to (-.5,.9) to (-.5,1.3) to (-.9,1.3) to (-.9,.9); 
\draw  (.7,.7) to (.7,.9);
\node at (.7,.5) {$\scriptstyle h_{(1)}$};
\draw (1,.7) to (.4,.7) to (.4,.3) to (1,.3) to (1,.7); 
\draw  (-.7,.7) to (-.7,.9);
\node at (-.7,.5) {$\scriptstyle h_{(0)}$};
\draw (-1,.7) to (-.4,.7) to (-.4,.3) to (-1,.3) to (-1,.7); 
\draw  (.7,.3) to [out=-90,in=0] (0,-.2);
\draw  (-.7,.3) to [out=-90,in=180] (0,-.2);
\draw  (0,-.2) to [out=-90,in=180] (.75,-.7);
\draw  (1.5,-.2) to  [out=-90,in=0] (.75,-.7);
\node at (3,.8) {$=\quad\ \sum_{(h)}$};
\end{tikzpicture}
\quad
\begin{tikzpicture}
\draw  (0,1.8) to  [out=0,in=90] (.7,1.3) to (.7,.3);
\draw  (0,1.8) to  [out=180,in=90] (-.7,1.3);
\node at (-.7,1.1) {$\scriptstyle a$};
\draw (-.9,.9) to (-.5,.9) to (-.5,1.3) to (-.9,1.3) to (-.9,.9); 
\draw  (-.7,.7) to (-.7,.9);
\node at (-.7,.5) {$\scriptstyle h_{(0)}$};
\draw (-1,.7) to (-.4,.7) to (-.4,.3) to (-1,.3) to (-1,.7); 
\draw  (.7,.3) to [out=-90,in=0] (0,-.2);
\draw  (-.7,.3) to [out=-90,in=180] (0,-.2);
\draw  (2,1.8) to  [out=0,in=90] (2.7,1.3);
\draw  (2,1.8) to  [out=180,in=90] (1.3,1.3) to (1.3,.3);
\node at (2.7,1.1) {$\scriptstyle b$};
\draw (2.9,.9) to (2.5,.9) to (2.5,1.3) to (2.9,1.3) to (2.9,.9); 
\draw  (2.7,.7) to (2.7,.9);
\node at (2.7,.5) {$\scriptstyle h_{(1)}$};
\draw (3,.7) to (2.4,.7) to (2.4,.3) to (3,.3) to (3,.7); 
\draw  (2.7,.3) to [out=-90,in=0] (2,-.2);
\draw  (1.3,.3) to [out=-90,in=180] (2,-.2);
\draw  (0,-.2) to [out=-90,in=180] (1,-.7);
\draw  (2,-.2) to  [out=-90,in=0] (1,-.7);
\end{tikzpicture}
$$
}}
\newline
coincides with (by Remark \ref{hap})
\newline
\scalebox{0.85}{\parbox{\linewidth}{%
$$
\begin{tikzpicture}
\node at (-1.8,.35) {$d(A)^{-1}\   \sum_{(h)} $};
\draw  (0,.9) to  [out=90,in=180] (.75,1.4);
\draw (-.2,.9) to (.2,.9) to (.2,.5) to (-.2,.5) to (-.2,.9); 
\node at (0,.7) {$\scriptstyle a$};
\draw  (0,.5) to (0,.2);
\draw (-.3,.2) to (.3,.2) to (.3,.-.2) to (-.3,-.2) to (-.3,.2); 
\node at (0,0) {$\scriptstyle h_{(0)}$};
\draw  (0,-.2) to [out=-90,in=180] (.75,-.7);
\draw   (.75,1.4) to [out=0,in=90] (1.5,.9);
\draw  (1.5,-.2) to  [out=-90,in=0] (.75,-.7);
\draw  (1.5,-.3) to (1.5,.9);
\end{tikzpicture}
\qquad
\begin{tikzpicture}
\draw  (0,.9) to  [out=90,in=180] (.75,1.4);
\draw (-.2,.9) to (.2,.9) to (.2,.5) to (-.2,.5) to (-.2,.9); 
\node at (0,.7) {$\scriptstyle b$};
\draw  (0,.5) to (0,.2);
\draw (-.3,.2) to (.3,.2) to (.3,.-.2) to (-.3,-.2) to (-.3,.2); 
\node at (0,0) {$\scriptstyle h_{(1)}$};
\draw  (0,-.2) to [out=-90,in=180] (.75,-.7);
\draw   (.75,1.4) to [out=0,in=90] (1.5,.9);
\draw  (1.5,-.2) to  [out=-90,in=0] (.75,-.7);
\draw  (1.5,-.3) to (1.5,.9);
\end{tikzpicture}
\quad.
$$
}}
\end{proof}

This implies that \'etale algebras have no ``quantum symmetries'' in the sense of 
Etingof-Walton \cite{ew}.
\begin{corollary}
Let $\C$ be a braided tensor category, and $A\in \C$ be an \'etale algebra. 
Assume that $H$ is a Hopf algebra faithfully acting on $A$.
Then $H$ is the group algebra $k[G]$ for some subgroup $G\subset \Aut_{\alg}(A)$.
\end{corollary}
\begin{proof}
Proposition \ref{prop:Qquotient} and \ref{prop:Qcomm} imply that $H^\ast$ is a quotient of 
the commutative algebra $Q(A)$ and therefore commutative.
Thus $H$ is co-commutative, which implies that $H=k[G]$ for some finite group $G$. 
Under this identification, restricting the action of $k[G]$ to  $G$ 
we obtain a group action of $G$ on $A$, and since the action is faithful, an embedding $G\to \Aut_{\alg}(A)$.
\end{proof}

\section{Maximally symmetric \'etale algebras in pseudo-unitary braided fusion categories}\label{gal}

Let now $\C$ be a braided pseudo-unitary fusion category.
Denote by $\Irr(\C)$ the set of (isomorphism classes of) simple objects of $\C$ and let $d_\C(X)$ be the (pseudo-unitary) dimension of $X \in  \C$.

\begin{lemma}\label{ml}
Let $A$ be an \'etale algebra.
Then
$$\dim(\C(X,A))\ \leq\ d(X)\ .$$
\end{lemma}
\begin{proof}
This follows from the chain of (in)equalities
$$\dim(\C(X,A)) = \dim(\C(I,X^*\otimes A)) = \dim(\C_A(A,X^*\otimes A)) \leq$$
$$\leq d_{\C_A}(X^*\otimes A) = d_\C(X^*)  = d_\C(X)\ .$$
The inequality is just the fact that multiplicity $\dim(\D(I,Y))$ of the unit object in another object $Y$ of a fusion category $\D(=\C_A)$ is bounded from above by the dimension $d_\D(Y)$.
\end{proof}
It follows from Lemma \ref{ml} that $\dim(\C(A,A)) \leq d(A)$ for any \'etale algebra $A$.
\newline
We call an \'etale algebra $A$ {\em maximally symmetric} if the above bound is saturated, i.e.\ if 
$$\dim(\C(A,A))\ =\ d(A)\ .$$

\begin{lemma}\label{al}
Let $A$ be a maximally symmetric \'etale algebra. Then $\dim(\C(X,A)) = d(X)$ for any $X\in \Irr(\C)$ such that $\C(X,A)\not=0$.
\end{lemma}
\begin{proof}
Write $A = \bigoplus_{X\in \Irr(\C)}\C(X,A)\otimes X$. Then 
$$\sum_{X\in \Irr(\C)}\dim(\C(X,A))^2 = \dim(\C(A,A)) = d(A) = \sum_{X\in \Irr(\C)}\dim(\C(X,A))d(X)\ .$$
Then Lemma \ref{ml} implies the desired.
\end{proof}

Denote by  $\Rep(G)$ the Tannakian (symmetric) category of finite dimensional representations of a group $G$. 
\begin{theorem}\label{mt}
Let $A$ be a maximally symmetric \'etale algebra. Let $G = \Aut_{\alg}(A)$ be the group of its algebra automorphisms. Then there is a full braided embedding $\Rep(G)\subset \C$ such that $A$ is isomorphic to the function algebra $k(G)$ considered as an algebra in $\Rep(G)\subset \C$. 
\end{theorem}
\begin{proof}
Consider the full subcategory $\E\subset\C$ additively generated by $X\in \Irr(\C)$ such that $\C(X,A)\not=0$. 
In other words, $\E$ is the support category of $A$ in $\C$, i.e.\ the smallest full subcategory containing $A$.
By Lemma \ref{al} the subcategory $\E$ coincides with the full subcategory 
$$
	\E = \{X\in \C ~|\ \ev_X\colon \C_A(A,X\otimes A)\otimes A\stackrel{\simeq}{\longrightarrow} X\otimes A\}
$$
of those $X$ for which the canonical (evaluation) morphism of right $A$-modules $\ev_X\colon \C_A(A,X\otimes A)\otimes A\to X\otimes A$ is an isomorphism.
The following property of evaluation maps shows that $\E$ is a (braided) tensor subcategory. 
For any $X, Y\in \C$ the diagram 
$$
\xymatrix{
(\C_A(A,X\otimes A)\otimes A)\otimes_A(\C_A(A,Y\otimes A)\otimes A) \ar[rr]^(.63){\ev_X\otimes_A\ev_Y} \ar@{=}[dd] && (X\otimes A)\otimes_A(Y\otimes A) \ar@{=}[d] \\
&& X\otimes Y\otimes A\\
\C_A(A,X\otimes A)\otimes \C_A(A,Y\otimes A)\otimes A \ar[rr] && \C_A(A,X\otimes Y\otimes A)\otimes A \ar[u]_{\ev_{X\otimes Y}}
}
$$
commutes.
Since the tensor product of morphisms $\C_A(A,X\otimes A)\otimes \C_A(A,Y\otimes A) \to \C_A(A,X\otimes Y\otimes A)$ is an injective linear map, the above diagram implies that the evaluation $\ev_{X\otimes Y}$ must be an isomorphism if $e_X$ and $e_Y$ are.
\newline
Now define a functor $F\colon\E\to\Vect$ by $F(X) = \C_A(A,X\otimes A) = \C(X^*,A)$. This functor is clearly faithful (since $A$ is self-dual as an object of $\C$).
It is also tensor, with the tensor structure
$$F(X)\otimes F(Y) = \C_A(A,X\otimes A)\otimes \C_A(A,Y\otimes A) \to \C_A(A,X\otimes Y\otimes A) = F(X\otimes Y)$$
given again by the tensor product of morphisms in $\C_A$.
Moreover this functor is braided, which makes the category $\E$ symmetric.
The Tannaka-Krein reconstruction with respect to $F$ provides the equivalence $\Rep(G)\to \E$, with $G = \Aut_\otimes(F)$.
Finally, note that the group homomorphism $\Aut_{\alg}(A)\to \Aut_\otimes(F)$, sending $g$ to the automorphism $\C(X^*,g)$ of $\C(X^*,A)$ is an isomorphism.
\end{proof}

\begin{corollary}
An  \'etale algebra is Galois if and only if it is maximally symmetric.
\end{corollary}
\begin{proof}
Let $A\in\C$ be a Galois \'etale algebra and let $G=\Aut_{\alg}(A)$. The functor $F\colon \Rep(G)\to\C$ defined by $F(U) = U\otimes_{k[G]}A$ is fully faithful. Indeed, 
$$\C(U\otimes_{k[G]}A,V\otimes_{k[G]}A) \simeq \Rep(G)(U,\C(A,V\otimes_{k[G]}A)) \simeq $$
$$\simeq \Rep(G)(U,V\otimes_{k[G]}\C(A,A))  \simeq \Rep(G)(U,V\otimes_{k[G]}k[G]) \simeq \Rep(G)(U,V)\ .$$
Since $F(k(G)) \simeq A$ and since $k(G)$ is a maximally symmetric algebra in $\Rep(G)$, we get that $A$ is maximally symmetric. 
\newline
Conversely since $k(G)$ is a Galois algebra in $\Rep(G)$, the embedding $\Rep(G)\subset \C$ associated with a maximally symmetric algebra $A\in\C$ makes $A$ Galois.
\end{proof}

For a simple complex Lie algebra $\g$ and a level $k\in\bZ_{>0}$ denote by $\g_k$ the vertex operator algebra structure on the simple vacuum module over the affinisation $\widehat\g$ at level $k$ \cite{fz}. Denote the tensor product of $\g_k$ and $\h_l$ by $\g_k\h_l$. 
Denote by $\C(\g,k)$ the ribbon fusion category of modules of $\g_k$, i.e. the category of representations of $\widehat\g$ at level $k$.
\begin{example}
Consider an embedding of vertex operator algebras (i.e. a conformal embedding) of the form $\sllie(3)_3\sllie(3)_3\subset\e_{8,1}$
(see e.g. \cite{bb,sw}). 
A conformal embedding $\sllie(3)_3\sllie(3)_3\subset\e_{8,1}$ corresponds to an \'etale algebra $L\in \C(\sllie(3),3)^{\boxtimes 2}$. 
Denote by $I=X_1,X_2,X_3$ the 1-dimensional objects of $\C(\sllie(3),3)$ and by $X$ the only 3-dimensional object of $\C(\sllie(3),3)$. It is not hard to see that the decomposition of $L$ into simple objects is
\begin{equation}\label{uol}L = (X_1\oplus X_2\oplus X_3)\boxtimes(X_1\oplus X_2\oplus X_3) \ \oplus\ 3(X\boxtimes X)\ .\end{equation}
The algebra $L$ is not maximally symmetric and hence is not Galois.
\newline
The algebra $L$ has an \'etale subalgebra $B = X_1\boxtimes X_1\oplus X_2\boxtimes X_2\oplus X_3\boxtimes X_3$, giving rise to a simple current extension 
$\sllie(3)_3\sllie(3)_3\subset \widetilde{\sllie(3)_3\sllie(3)_3}$ of index 3. 
The category of modules over the vertex operator algebra $\widetilde{\sllie(3)_3\sllie(3)_3}$ coincides with the category $(\C(\sllie(3),3)^{\boxtimes 2})^{\loc}_B$ of local $B$-modules in $\C(\sllie(3),3)^{\boxtimes 2}$ (see \cite{hkl,ckm}). 
The extension $B\subset L$ gives rise to an \'etale algebra $\overline L$ in $(\C(\sllie(3),3)^{\boxtimes 2})^{\loc}_B$.
 Its not hard to show that there are 4 such algebras with decompositions into simple objects 
\begin{equation}\label{dea}\overline L_i = Y_1\oplus Y_2\oplus Y_3\oplus 3Z_i,\qquad \overline L = Y_1\oplus Y_2\oplus Y_3\oplus Z_1\oplus Z_2\oplus Z_3\ ,\end{equation}
where $Y_i$ are the 1-dimensional simple $B$-modules induced from $X_i\boxtimes X_i$ and $Z_i$ are the 3-dimensional simple submodules of the $B$-modules induced from $X\boxtimes X$. 
Note that these algebras considered as algebra in $\C(\sllie(3),3)^{\boxtimes 2}$ all have the same underlying object \eqref{uol}.
\newline
The first type algebra $\overline L_i$ is maximally symmetric and is Galois in $(\C(\sllie(3),3)^{\boxtimes 2})^{\loc}_B$.
The decomposition \eqref{dea} shows that the automorphism group $\Aut_{\alg}(\overline L_i)$ is the alternating group $A_4$.
Thus we have shown that the simple current extension 
$\widetilde{\sllie(3)_3\sllie(3)_3}$ is the vertex operator subalgebra of invariants under an $A_4$-action:
$$(\e_{8,1})^{A_4} = \widetilde{\sllie(3)_3\sllie(3)_3}\ .$$
\end{example}


\end{document}